\documentclass[notitlepage,12pt]{amsart}%
\usepackage{geometry}
\usepackage{natbib}
\usepackage[singlespacing]{setspace}
\usepackage{graphicx}
\usepackage{amsmath}
\usepackage{amsfonts}
\usepackage{amssymb}%
\newtheorem{theorem}{Theorem}

\begin{document}

\title{Incentive-Compatible Elicitation of Quantiles}
\author{Nicholas M. Kiefer }
\thanks{Cornell University, Departments of Economics and
	Statistical Science, 490 Uris Hall, Ithaca, NY 14853-7601, US,
	email:nicholas.kiefer@cornell.edu, and CREATES, funded
	by the Danish Science Foundation, University of Aarhus,Denmark.
} 

\date{\today}
\begin{abstract}
Incorporation of expert information in inference or decision settings is often
important, especially in cases where data are unavailable, costly or
unreliable. One approach is to elicit prior quantiles from an expert and then
to fit these to a statistical distribution and proceed according to Bayes
rule. Quantiles are often thought to be easier to elicit than moments. An
incentive-compatible elicitation method using an external randomization is available. Such a mechanism will encourage the expert to exert the care necessary to report accurate information. 
A second application might be called posterior elicitation. Here an analysis has been done and the results must be reported to a decision maker. For a variety of reasons (possibly including the reward system in the corporate hierarchy) the modeler might need the right incentive system to report results accurately. Again, eliciting posterior quantiles can be done with an incentive compatible mechanism.

\textbf{MSC 2000 Subject classification}, Primary: 62C10; secondary 91B06

\textbf{Keywords}: Bayesian inference, mechanism design, prior assessment

\end{abstract}
\maketitle
\section{Introduction}

Incorporation of prior information is important in any decision or inference
setting, whether it is done formally or informally. The formal Bayesian
approach encourages transparency in assumptions, clear thinking and coherence.
One example is risk management in financial institutions in which prudent
management requires understanding default probabilities for groups of
homogeneous assets. In the case of  new types of assets there may
not be enough data information to support a practical conventional estimator,
for example the frequency estimator in the case of binomial defaults. Zero is not an estimated default probability that is acceptable to regulators. This
issue has attracted regulatory and industry as well as academic attention, see
\cite{Kiefer2009a}, and references given there. \cite{Kiefer2010} proposes
eliciting prior quantiles for an expert's prior on the value of the default
probability for a particular group of assets. For example, the median can be
assessed by asking the expert at what value of a default rate $\theta$ would
he be equally surprised to see a realization above or below $\theta.$ These
quantiles (perhaps after feedback and revision) are assembled into a
distribution, either by fitting a specific functional form or as proposed by
\cite{Kiefer2010} fit to a smoothed maximum-entropy distribution. The idea is
to impose as little information as possible beyond that elicited from the
expert. This distribution is then used to process data information through
Bayes rule and the likelihood function. The latter is itself typically a
representation of a large number of probabilities in terms of a few
parameters, so the statistical treatment and the approximations involved are
the same in the prior and the likelihood.

Other examples of elicitation of quantiles and their use to form a prior
distribution are cited in \cite{O'Hagan2006} and include applications to drug
testing, sales of engines, the effect of nuclear waste on temperature, and
future earnings. A discussion of the statistical and psychological issues involved in elicitation is \citet{Garthwaite2005}. These issues are not reviewed here. The incentive compatibility question does not seem to have
been stressed.

A second application is eliciting accurate assessments from a modeling group within an organization. Here, the risks and rewards to proper reporting incentives are clear. In banks, risk modelers are required by  regulators to report to senior management, not to business line management. In other institutions the separation between incentives for modelers and those acting on the results might be less clear. Senior management may choose to put in place incentive compatible reporting mechanisms. Regulators may question the compensation plan for modelers.

Finally, reporting of financial forecasts and results to counterparties and regulators should be as accurate as possible. Here too, incentive-compatible mechanisms might play a role. 

The difficult part is the elicitation of the quantiles, which requires thought
and therefore some effort from the expert in our first example (the expert must have the incentive to provide this effort) and simply rewards for accuracy in the second. Since the quantiles can never be
observed, and therefore the assessment checked, there is an issue of providing
an incentive for the expert to provide the required thought. The problem of
eliciting probabilities for given sets is well-studied and a widely-used
approach is \textit{scoring}. The scoring method does not naturally extend to
the quantile assessment problem, as shown in Section 2. \cite{Savage1971}
reviews techniques for assessing probabilities and notes an interesting
interpretation of probabilities as prices. He notes that the device of outside
randomization, used by \cite{Marschak1964} to compel a true valuation for a
bid or asked price applies also to probability assessment. An ingenious recent
method due to \cite{Karni2009} introduces a second outside source of
randomness to eliminate possible effects of risk aversion. This method does
extend naturally to eliciting quantiles as shown in Section 2 with a different
development than that of Marschak or Karni. Section 3 concludes.

\section{Eliciting Probabilities and Quantiles}

Assume at the outset that the expert's information about the unknown quantity
$\theta$ is coherent, that is that it can be described by a probability
distribution. Classical discussions of the necessity of describing uncertainty
in terms of probability are \cite{Savage1954}, \cite{De1974}, and
\cite{Lindley1982a}. We do not review these well-known demonstrations and
simply assume that the expert's information is described in a probability
distribution with $cdf$ $F(\theta)$ and $pdf$ $f(\theta)$. We wish to elicit
the probability $\alpha$ quantile $q_{\alpha}$ with $F(q_{\alpha})=\alpha.$
Assume that $f(q_{\alpha})>0$ so that the quantile is well-defined. Assume for
convenience that $Supp(f)\subseteq\lbrack0,1].$ This is natural when the
uncertain quantity $\theta$ is a default probability; in other cases a
parameter transformation may be appropriate. The method of scoring for
eliciting the probability $\beta$ associated with a given set, say $[0,q]$,
rewards the expert with
\begin{equation}
r(\beta|q)=(I(\theta\in\lbrack0,q])g_{1}(\beta)+(1-I(\theta\in\lbrack
0,q])g_{0}(\beta)\label{probscore}%
\end{equation}
where $I()$ is the indicator function $g_{1}(\beta)$ is a nondecreasing
function, $g_{0}(\beta)$ is a nonincreasing function and $\beta$ is the
elicited probability, after the single realization of $\theta$ is seen. The
scoring rule is proper if the expectation $Er(\beta|q)=F(q)g_{1}%
(\beta)+(1-F(q))g_{0}(\beta)$ is maximized at $\beta=q_{\alpha}$ (some
applications minimize instead)$.$ The optimal choice of $\beta$ for a
risk-neutral expert rewarded by a proper scoring rule is $F(q).$ A classical
example of a proper scoring rule has $g_{1}(\beta)$ $=(1-\beta)^{2}$ and
$g_{0}(\beta)=-\beta^{2}$ in which case the first-order condition for
maximizing $Er(\beta|q)$ implies $F(q)/(1-F(q))=\beta/(1-\beta).$

Scoring rules for probabilities have been widely studied. A classical
application is the assessment of the quality of weather forecasts. Most of the
literature does not consider scoring for probability assessment, rather for
estimation when the score is an objective function or for measuring the
accuracy of assessments already arrived at, or for evaluating the fit of
statistical models. Thus, most of the literature does not consider risk
aversion and its effects on probability assessment, although the issue is well
known, see \cite{Savage1971} or for a general treatment that also considers
the effects of state preference (a very difficult problem) \cite{Kadane1988}.
\cite{Schervish1989} provides a characterization of the class of proper
scoring rules for probabilities. It is clear from (\ref{probscore}) that
scoring rules need not be symmetric. \cite{Winkler1994} argues that in
forecasting a short-horizon weather event the maximum score should not occur
at probability 0.5 as with symmetric scores but at the long-run event
probability, reflecting maximum uncertainty. This argument may also be
relevant for assessing financial risks. Scoring the assessment of the
probability of a binary event can be readily extended to scoring the
assessment of a full probability distribution. A simple device is to calculate
the score for a distribution assessment by choosing an interval randomly and
scoring the probability implied by the assessed distribution using a scoring
rule for a binary event. \cite{Matheson1976} propose a score that integrates
over the random interval with a weighting function emphasizing more important
parts of the distribution being assessed.  Scoring rules for probabilities can
be related to information measures and utilities, see \cite{Jose2008}. Proper
scoring rules typically do not lead to accurate assessments in the presence of
risk aversion. 

With $u(x)$ the utility of a payoff $x$ the expected utility associated with
the scoring rule (\ref{probscore}) is
\[
Eu(r(\beta|q))=F(q)u(g_{1}(\beta))+(1-F(q))u(g_{0}(\beta))
\]
and the expected-utility maximizing choice of $\beta$ is not typically $F(q)$
when $u$ is nonlinear. In the quadratic case for example the first-order
condition implies $F(q)/(1-F(q))=u^{\prime}(-\beta^{2})\beta/(u^{\prime
}((1-\beta)^{2})(1-\beta)).$ \cite{Karni2009} provides a method for eliciting
probabilities that works in the presence of risk aversion. That method extends
to quantile elicitation. 

There is far less work on quantile scoring. Fixing $\beta$ and using
(\ref{probscore}) to elicit the quantile $q$ does not work even without risk
aversion as the optimal choice of $q$ is 0 for $g_{0}(\beta)\geq g_{1}(\beta)$
and 1 for $g_{0}(\beta)\leq g_{1}(\beta).$ \cite{Gneiting2007} note that
$w(q;\beta)=\beta s(q)+(s(\theta)-s(q))1\{\theta\leq q\}+h(\theta)$ for $s$
nondecreasing and $h$ arbitrary gives a proper scoring rule for quantiles.
They note that the characterization of the full class of proper scoring rules
for quantiles remains open. We give a simple proof for differentiable and
increasing $s.$ 

\begin{theorem}
$w(q|\beta)=\beta s(q)+(s(\theta)-s(q))1\{\theta\leq q\}+h(\theta)$ with s()
increasing and h() arbitrary is a proper scoring rule for the $\beta th$
quantile $q.$

\begin{proof}
$E_{\theta}w(q|\beta)=\beta s(q)+\int\nolimits_{0}^{q}(s(\theta)-s(q))dF+\int
\nolimits_{0}^{1}h(\theta)dF.$ The first order condition is $\beta s^{\prime
}(q)-F(q)s^{\prime}(q)=0,$ hence the optimal q satisfies $F(q)=\beta$. 
\end{proof}
\end{theorem}

A widely used rule in econometrics (for estimation and goodness of fit
assessment, not for eliciation) is $w(q|\beta)=(\theta-q)(1\{\theta\leq
q\}-\beta)$. This scoring rule, with $s(q)=q$, is behind most work in quantile estimation, see
\cite{Koenker1978} and \cite{Koenker1999}. With risk aversion, proper scoring
rules for quantiles need not elicit the true quantiles. 

\section{A New Method for Quantile Elicitation}

The new method based on outside randomness works as follows. Suppose the
elicitor wishes $q_{a}$. The elicitor has access to a genie which generates a
random variable $\xi$ from the uniform distribution on $[0,1]$ and
independently a random variable $d\in\{0,1\}$ from the Bernoulli with
probability $\alpha.$  The expert
supplies a value $q$ for the $\alpha-th$ quantile. Nature supplies one realization of $\theta.$ The expert receives a
reward equal to $rI(\theta\in\lbrack0,\xi])$ if $q<\xi$ and $rd$ if $q\geq
\xi.$  The expert's utility of a
payoff $x$ is $u(x).$

It is the random variable $ \xi $ which leads to true revelation in the risk neutral case and the Bernoulli $d$ which allows for risk aversion. To develop intuition, consider the method  with $ d $ replaced by its expectation $ \alpha $.Then the reward from 4) upon observing $ \theta $ is 
\begin{equation*}
v(q|\xi,\theta)=u(r)I(\theta\in [0,\xi])I(q<\xi)+u(r\alpha)I(q\geq\xi))
\end{equation*}
(normalizing $u(0)=0$). Marginalizing wrt $ \theta $ and then $ \xi $ leads to 
\begin{equation*}
v(q)=u(r)\int_{q}^{1}F(t)dt +u(r\alpha)q
\end{equation*}
and the FOC $u(r)(-F(q))+u(r\alpha)$ leads to $q=q_\alpha$ when utility is linear but not otherwise.

Now suppose the rv $d$ is supplied along with $ \xi $. Then
\begin{equation*}
v(q|\xi,d)=u(r)F(\xi)I(q<\xi)+u(rd)I(q\ge \xi)
\end{equation*}
and marginalizing wrt $d$ gives
\begin{equation*}
v(q|\xi)=u(r)F(\xi)I(q<\xi)+u(r)I(q\ge \xi)\alpha. 
\end{equation*}
The FOC implies $F(q)=\alpha$, so revelation is optimal with risk aversion. The random variable $d$  supplied by the genie essentially moves the expectation through the utility function. In short:
\begin{theorem}
With the reward described above the optimal policy for the expert is to report
the true quantile.

\begin{proof}
Consider the expected utility to the expert of supplying $q.$ First, suppose
the genie supplies $\xi$ to the expert. Marginalizing wrt $\theta$ and $d$ gives the expected utility 
\[
v(q|\xi)=u(r)(F(\xi)I(q<\xi)+F(q_{\alpha})I(q\geq\xi))
\]
piecewise constant with a break at $\xi.$ Marginalizing with respect to the
uniform random variable $\xi$ yields the unconditional expected utility
function%
\[
v(q)=u(r)(\int_{q}^{1}F(t)dt+F(q_{\alpha})q)
\]
The first-order condition is $v^{\prime}(q)=u(r)(-F(q)+F(q_{\alpha}))=0$ and
the function is concave, so the optimal policy for the expert is to report the
true quantile.
\end{proof}
\end{theorem}

An alternative proof from a decision-theoretic point of view and using
lotteries can be given. This proof uses preferences over lotteries but does
not require the full expected utility framework. Let $(x,p)\in R\times
\lbrack0,1]$ denote the lottery that pays $\$x$ with probability $p$ and \$0
with probability $(1-p)$. The expert payoff is $(r,F(\xi))$ if $\xi>q$ and
$(r,\alpha)$ if $\xi\leq q$. Consider the report $q>q_{\alpha}$. If $\xi>q$
then the expert's payoff is $(r,F(\xi))$ whether he reports $q$ or $q_{\alpha
}$. If $\xi\leq q_{\alpha}$ then the expert's payoff is $(r,\alpha)$ whether
he reports $q$ or $q_{\alpha}$. If $q_{\alpha}<\xi$ $<q$ then the expert's
payoff is $(r,\alpha)$. However, had he reported $q_{\alpha}$ instead, his
reward would have been $(r,F(\xi))$ . But $F(\xi)>\alpha$ hence $(r,F(\xi))$
first-order stochastically dominates $(r,\alpha)$ so the expert cannot win and
may lose as a result of reporting $q>q_{\alpha}.$Similarly, reporting
$q<q_{\alpha}$ is dominated.

\section{Conclusion}

The classical elicitation problem concerns eliciting probabilities for given
events. This paper studies the complementary problem of eliciting events for
given probabilities. This is the problem involved in obtaining prior
quantiles. Although the reward $r$ does not affect the optimality condition,
it is clear that the actual effort expended by the expert will depend on the
value of the reward (precisely, on its utility). Perhaps some part of a bonus
could be tied into the probability assessment. Interesting open questions
include: Which quantiles and how many should be assessed? How much accuracy
can be expected in a quantile assessment? Can experts be trained to improve
their assessments? How can prior quantiles be assessed from a group of experts?

\bibliographystyle{econometrica}
\bibliography{C:/AAADocuments/Utilities/banking,C:/AAADocuments/Utilities/Bayes,C:/AAADocuments/Utilities/Kiefer2}

\end{document}